\documentclass[11pt]{article}
\usepackage{amssymb}

\newtheorem{theorem}{Theorem}
\newtheorem{proposition}{Proposition}
\newtheorem{lemma}{Lemma}
\newtheorem{corollary}{Corollary}

\newcommand{\Z}{{\mathbb Z}}

\newcommand{\B}{{\cal B}}
\newcommand{\R}{{\mathbb R}}
\newcommand{\C}{{\mathbb C}}
\newcommand{\EE}{{\mathbb E}}
\newcommand{\F}{{\cal F}}
\newcommand{\LL}{\Lambda_3}
\newcommand{\PP}{\mathbb{P}}
\newcommand{\SIGMA}{\raisebox{-0.4ex}{\mbox{\Large $\Sigma$}}}

\title{The Minimal Number of Three-Term Arithmetic Progressions 
Modulo a Prime Converges to a Limit}

\author{Ernie Croot \thanks{Supported by an NSF grant}}

\begin{document}

\maketitle

\begin{abstract} How few three-term arithmetic progressions can a
subset $S \subseteq \Z_N := \Z/N\Z$ have if $|S| \geq \upsilon N$?
(that is, $S$ has density at least $\upsilon$).  
Varnavides \cite{varnavides} 
showed that this number of arithmetic-progressions is at
least $c(\upsilon)N^2$ for sufficiently large integers $N$; 
and, it is well-known that determining good lower bounds for 
$c(\upsilon)> 0$  is at the same level of depth as Erd\" os's famous
conjecture about whether a subset $T$ of the naturals where
$\sum_{n \in T} 1/n$ diverges, has a $k$-term arithmetic progression
for $k=3$ (that is, a three-term arithmetic progression).   

The author answers a question of B. Green \cite{AIM} 
about how this minimial number of progressions oscillates
for a fixed density $\upsilon$ as $N$ runs through the primes, and
as $N$ runs through the odd positive integers.
\end{abstract}

\section{Introduction}

Given an integer $N \geq 2$ and a mapping $f : \Z_N \to \C$ 
define
\begin{eqnarray}
\LL(f)\ =\ \LL(f;N)\ &:=&\ \EE_{n,d \in \Z_N}(f(n)f(n+d)f(n+2d))
\nonumber \\
&=&\ {1 \over N^2} \sum_{n,d \in \Z_N} f(n)f(n+d)f(n+2d), \nonumber
\end{eqnarray}
where $\EE$ is the expectation operator, defined for a function
$g : \Z_N \to \C$ to be 
$$
\EE(g)\ =\ \EE_n(g)\ :=\ {1 \over N} \sum_{n \in \Z_N} g(n).
$$
 
If $S \subseteq \Z_N$, and if we identify $S$ with its indicator
function $S(n)$, which is $0$ if $n \not \in S$ and is $1$ if $n \in S$,
then $\LL(S)$ is a normalized count of the number of three-term 
arithmetic progressions $a,a+d,a+2d$ 
in the set $S$, including trivial progressions $a,a,a$.   
\bigskip

Given $\upsilon \in (0,1]$, consider the family $\F(\upsilon)$ 
of all functions 
$$
f : \Z_N \to [0,1],\ {\rm such\ that\ }\EE(f) \geq \upsilon.
$$
Then, define
$$
\rho(\upsilon,N)\ :=\ \min_{f \in \F(\upsilon)} \LL(f).
$$  
From an old result of Varnavides \cite{varnavides} we know that 
$$
\LL(f)\ \geq\ c(\upsilon)\ >\ 0,
$$
where $c(\upsilon)$ does not depend on $N$.  A natural and interesting question 
(posed by B. Green \cite{AIM}) is to determine whether for fixed $\upsilon$
$$
\lim_{p \to \infty \atop p\ {\rm prime}} \rho(\upsilon,p)\ {\rm exists}?
$$

In this paper we answer this question in the affirmative:
\footnote{The harder, and more interesting 
question, also asked by B. Green, which we do not answer in this
paper, is to give a simple formula for this limit.}

\begin{theorem}  \label{main_theorem} For a fixed $\upsilon \in (0,1]$ we have
$$
\lim_{p \to \infty \atop p\ {\rm prime}} \rho(\upsilon,p)\ {\rm exists}.
$$
\end{theorem}

Call the limit in this theorem $\rho(\upsilon)$.  Then, this theorem has the 
following immediate corollary:

\begin{corollary} For a fixed $\upsilon \in (0,1]$, let $S$ be any subset of $\Z_N$
such that $\LL(S)$ is minimal subject to the constraint $|S| \geq \upsilon N$.
Let $\rho_2(\upsilon, N) = \LL(S)$.  Then, 
$$
\lim_{p \to \infty \atop p\ {\rm prime}} \rho_2(\upsilon,p)\ =\ \rho(\upsilon).
$$
\end{corollary}

Given Theorem \ref{main_theorem}, the proof of the corollary is standard, 
and just amounts to applying a functions-to-sets
lemma, which works as follows:  Given $f : \Z_N \to [0,1]$, we let $S_0$ be a random subset of 
$\Z_N$ where $\PP(s \in S_0) = f(s)$.  It is then easy to show that with probability $1-o(1)$,
$$
\EE(S_0)\ \sim\ \EE(f),\ {\rm and\ } \LL(S_0)\ \sim\ \LL(f).
$$
So, there will exist a set $S_1$ with these two properties (an instantiation
of the random set $S_0$).  Then, by adding only a 
small number of elements to  $S_1$ as needed, we will have a set $S$ satisfying
$$
|S|\ \geq\ \upsilon N,\ {\rm and\ } \LL(S)\ \sim\ \LL(f).
$$
\bigskip

We will also prove the following:

\begin{theorem} \label{second_theorem} For $\upsilon = 2/3$ we have that 
$$
\lim_{N \to \infty \atop N\ {\rm odd}} \rho(\upsilon, N)\ {\rm does\ not\ exist},
$$
where here we consider all odd $N$, not just primes. 
\end{theorem}

Thus, in our proof of Theorem \ref{main_theorem} we will make special use of the fact
that our moduli are prime.
\bigskip

\section{Basic Notation on Fourier Analysis}

Given an integer $N \geq 2$ (not necessarily prime), 
and a function $f : {\mathbb Z}_N \to {\mathbb C}$, we
define the Fourier transform
$$
\hat f(a)\ =\ \sum_{n \in \Z_N} f(n) e^{2\pi i a n / N}.
$$
Thus, the Fourier transform of an indicator function $C(n)$ for a set $C \subseteq \Z_N$
is:
$$
\hat C(a)\ =\ \sum_{n=0}^{N-1} C(n) e^{2\pi i an/N}\ =\ \sum_{n \in C} 
e^{2\pi i an/N}.
$$

Throughout the paper, when working with Fourier transforms, 
we will use a slightly compressed form of summation notation,
by introducing the sigma operator, defined by 
$$
\SIGMA_n\  f(n)\ =\ \sum_{n \in \Z_N} f(n).
$$
\bigskip

We also define define the norms
$$
||f||_t\ =\ (\EE |f(n)|^t )^{1/t},
$$
which is the usual $t$-norm where we take our measure to be
the uniform measure on $\Z_N$.

With our definition of norms, H\"older's inequality takes the form
$$
||f_1 f_2 \cdots f_n||_b\ \leq\ ||f_1||_{b_1} ||f_2||_{b_2} 
\cdots ||f_n||_{b_n}, \ {\rm if\ } {1 \over b}\ =\ {1 \over b_1} + 
\cdots + {1 \over b_n},
$$
although we will ever only need this for the product of two functions, and where the
$a_i$ and $b_i$ are $1$ or $2$ (i.e. Cauchy-Schwarz).

In our proofs we will make use of Parseval's identity, which says that
$$
||\hat f||_2^2\ =\ N ||f||_2^2
$$
This implies that 
$$
||\hat C||_2^2\ =\ N |C|.
$$
We will also use Fourier inversion, which says
$$
f(n)\ =\ N^{-1} \SIGMA_a e^{-2\pi a n/N} \hat f(a).
$$
Another basic fact we will use is that 
$$
\LL(f)\ =\ N^{-3} \SIGMA_a\ \hat f(a)^2 \hat f(-2a).
$$

\section{Key Lemmas}

Here we list some key lemmas we will need in the course of our proof
of Theorems \ref{main_theorem} and \ref{second_theorem}.

\begin{lemma} \label{large_values}
Suppose $h : \Z_N \to [0,1]$, and let ${\cal C}$ denote the set of all values 
$a \in \Z_N$ for which 
$$
|\hat h(a)|\ \geq\ \beta \hat h(0).
$$
Then,
$$
|{\cal C}|\ \leq\ (\beta \hat h(0))^{-2}N^2.
$$
\end{lemma}
\bigskip

\noindent {\bf Proof of the Lemma.}  This is an easy consequence of Parseval:
$$
|{\cal C}| (\beta \hat h(0))^2\ \leq\ N||\hat h||_2^2\ =\ N^2 ||h||_2^2\ 
\leq N^2.
\ \ \ \blacksquare
$$
\bigskip

\begin{lemma} \label{same_fourier} 
Suppose that $f, g : \Z_N \to [-2,2]$ have the 
property
$$
||\hat f - \hat g||_\infty\ <\ \beta N.
$$
Then, 
$$
|\LL(f) - \LL(g)|\ <\ 12 \beta.
$$
\end{lemma} 
\bigskip

\noindent {\bf Proof of the Lemma.}  The proof is an exercise in multiple uses of
Cauchy-Schwarz (or H\"older's inequality) and Parseval.  

First, let $\delta(a) = \hat f(a) - \hat g(a)$.  We have that 
\begin{eqnarray}
\LL(f)\ &=&\ N^{-3} \SIGMA_a \hat f(a)^2 (\hat g(-2a) + \delta(-2a)) \nonumber \\
&=&\ N^{-3} \SIGMA_a \hat f(a)^2 \hat g(-2a)\ +\ E_1, \nonumber
\end{eqnarray}
where by Parseval's identity we have that the error $E_1$ satisfies
$$
|E_1|\ \leq\  N^{-2} ||\delta||_\infty ||\hat f||_2^2\ =\ N^{-1} ||\delta||_\infty ||f||_2^2
\ <\ 4\beta.
$$

Next, we have that 
\begin{eqnarray}
N^{-3} \SIGMA_a \hat f(a)^2 \hat g(-2a)\ &=&\ N^{-3}  \SIGMA_a \hat f(a) 
(\hat g(a) + \delta(a)) \hat g(-2a) \nonumber \\
&=&\ N^{-3} \SIGMA_a \hat f(a) \hat g(a) \hat g(-2a)\ +\ E_2,
\nonumber
\end{eqnarray}
where by Parseval again, along with Cauchy-Schwarz (or H\"older's inequality), 
we have that the error $E_2$ satisfies

$$
|E_2|\ \leq\ N^{-2} ||\hat f(a) \hat g(-2a)||_1 ||\delta||_\infty 
\ <\ \beta N^{-1} ||\hat f||_2 ||\hat g||_2\ \leq\ 4\beta. 
$$

Finally,
\begin{eqnarray}
N^{-3} \SIGMA_a \hat f(a) \hat g(a) \hat g(-2a)\ &=&\ 
N^{-3} \SIGMA_a (\hat g(a) + \delta(a)) \hat g(a) \hat g(-2a) \nonumber \\
&=&\ \LL(g)\ +\ E_3, \nonumber
\end{eqnarray}
where by Parseval again, along with Cauchy-Schwarz (H\"older), 
we have that the error $E_3$ satisfies

$$
|E_3|\ \leq\ N^{-2} ||\delta||_\infty ||\hat g(a) \hat g(-2a)||_1\ <\ \beta N^{-1} 
||\hat g||_2^2\ =\ \beta ||g||_2^2\ \leq\ 4\beta.
$$

Thus, we deduce
$$
|\LL(f) - \LL(g)|\ <\ 12 \beta.\ \ \ \blacksquare
$$
\bigskip

The following Lemma and the Proposition after it make use of ideas
similar to the ``granularization'' methods from \cite{green1} and
\cite{green2}.

\begin{lemma} \label{flatten_lemma} 
For every $t \geq 1$, $0 < \epsilon < 1 $, 
the following holds for all primes $p$ sufficiently large: 
Given any set of residues $\{b_1,...,b_t\} \subset \Z_p$, there
exists a weight function $\mu : \Z_p \to [0,1]$ such that

$\bullet$ $\hat \mu(0)\ =\ 1$ (in other words, $\EE(\mu) = p^{-1}$);

$\bullet$ $|\hat \mu(b_i) - 1| < \epsilon^2$, for all $i=1,2,...,t$; and,

$\bullet$  $|| \hat \mu ||_1\ \leq\ p^{-1} (6\epsilon^{-1})^t$.
\end{lemma}

\noindent {\bf Proof.}  We begin by defining the functions
$y_1,...,y_t : \Z_p \to [0,1]$ by defining their Fourier transforms:  
Let $c_i \equiv b_i^{-1} \pmod{p}$, $L = \lfloor \epsilon p / 10 \rfloor$,
and define
$$
\hat y_i(a)\ =\ (2L+1)^{-1} \left ( \SIGMA_{|j| \leq L}  
e^{2\pi i a c_i j/p} \right )^2\ \in\ {\mathbb R}_{\geq 0}.
$$
It is obvious that $0 \leq y_i(n) \leq 1$, and $y_i(0) = 1$.  Also note that 
\begin{equation} \label{yi_fact}
y_i(n) \neq 0\ \ {\rm implies\ \ } b_i n \equiv j  \pmod{p},\ {\rm where\ }
|j| \leq 2L.
\end{equation}

Now we let $v(n) = y_1(n) y_2(n) \cdots y_t(n)$.  Then,
\begin{eqnarray} \label{v_transform}
\hat v(a)\ &=&\ p^{-t+1}(\hat y_1 * \hat y_2 * \cdots * \hat y_t)(a) \nonumber \\
&=&\ p^{-t+1} \SIGMA_{r_1+\cdots + r_t \equiv a}\ 
\hat y_1(r_1) \hat y_2(r_2) \cdots \hat y_t(r_t).
\end{eqnarray}
Now, as all the terms in the sum are non-negative reals we deduce that for $p$
sufficiently large,
\begin{eqnarray} \label{v0_bound}
p\ >\ \hat v(0)\ \geq\ p^{-t+1} \hat y_1(0)\cdots \hat y_t(0)\ &=&\ 
p^{-t+1} (2L+1)^t\nonumber \\ 
&>&\ (\epsilon/6)^t p.
\end{eqnarray}

We now let $\mu(a)$ be the weight
whose Fourier transform is defined by
\begin{equation} \label{w_def}
\hat \mu(a)\ =\ \hat v(0)^{-1} \hat v(a).
\end{equation}
Clearly, $\mu(a)$ satisfies conclusion 1 of the lemma.  

Consider now the value $\hat \mu(b_i)$.  As $\mu(n) \neq 0$ implies $y_i(n) \neq 0$,
from (\ref{yi_fact}) we deduce that if $\mu(n) \neq 0$, then for some $|j| \leq 2L$,
$$
{\rm Re}(e^{2\pi i b_i n/p})\ =\ {\rm Re}(e^{2\pi i j / p})\ =\ 
\cos(2\pi j / p)\ \geq\ 
1 - {1 \over 2} (2\pi \epsilon/5)^2\ >\ 1 - \epsilon^2. 
$$
So, since $\hat \mu(b_i)$ is real, we deduce that 
$$
\hat \mu(b_i)\ =\ \hat v(0)^{-1} \SIGMA_n v(n) 
e^{2\pi i b_i n/p}\ >\ 1-\epsilon^2.
$$
So, our weight $\mu(n)$ satisfies the second conclusion of our Lemma.  

Now, then, from (\ref{v_transform}), (\ref{w_def}), and
(\ref{v0_bound}) we have that
\begin{eqnarray} 
||\hat u||_1\ &=&\ p^{-t}\hat v(0)^{-1}  
\SIGMA_a\ \SIGMA_{r_1+\cdots + r_t \equiv a}\ 
\hat y_1(r_1) \hat y_2(r_2)\cdots \hat y_t(r_t) \nonumber \\
&=&\ p^{-t} v(0)^{-1} \prod_{i=1}^t \SIGMA_r \hat y_i(r) 
\nonumber \\
&=&\ \hat v(0)^{-1} y_1(0)y_2(0)\cdots y_t(0) \nonumber \\
&=&\ \hat v(0)^{-1} \nonumber \\
&<&\ p^{-1}(6\epsilon^{-1})^t.\ \ \ \ \blacksquare \nonumber
\end{eqnarray}
\bigskip

Next we have the following Proposition, which is an extended corollary
of Lemmas \ref{same_fourier} and \ref{flatten_lemma}:

\begin{proposition} \label{measure_prop}  For every $\epsilon > 0$,
$p > p_0(\epsilon)$ prime, and every $f : \Z_p \to [0,1]$,
there exists a periodic function $g : \R \to \R$ with period $p$ satisfying:

$\bullet$ $\EE(g) = \EE(f)$ (Here when we compute the expectation of $g$ we
restrict to $g : \{0,...,p-1\} \to \R$, and treat it as a mapping
from $\Z_p$.)

$\bullet$ $g : \R \to [-2\epsilon, 1+2\epsilon]$.

$\bullet$  $\hat g$ has ``small'' (approximate) support, when treated
as a function from $\Z_p \to \R$.   That is, 
there is a set of residues $c_1,...,c_m \in \Z_p$, $m < m_0(\epsilon)$, satisfying
$$
g(n)\ =\ p^{-1} \SIGMA_{1 \leq i \leq m} e^{-2\pi i c_i n/p} \hat g(c_i).
$$

$\bullet$ The $c_i$ satisfy $|c_i| < p^{1-1/m}$.

$\bullet$ $|\LL(g) - \LL(f)| < 25\epsilon$. 
\end{proposition}
\bigskip

\noindent {\bf Proof of the Proposition.}  We will need to define a number of sets and
functions in order to begin the proof:  Define
$$
\B\ =\ \{a \in \Z_p\ :\ |\hat f(a)| > \epsilon \hat f(0)\},
$$
and let $t = |\B|$.  Define
$$
\B'\ =\ \{a \in \Z_p\ :\ |\hat f(-2a)|\ {\rm or\ } |\hat f(a)|\  >\ \epsilon (\epsilon/6)^t \hat f(0) \},
$$
and let $m = |\B'|$.   Note that $\B \subseteq \B'$ implies $t \leq m$.  
Lemma \ref{large_values} implies that $m < m_0(\epsilon)$,
where $m_0(\epsilon)$ depends only on $\epsilon$.

Let $\mu : \Z_p \to [0,1]$ be as in Lemma \ref{flatten_lemma} with parameter $\epsilon$
and with $\{b_1,...,b_t\} = \B$.

Let $1 \leq s \leq p-1$ be such that for every $b \in \B'$,
$$
b\ \equiv\ s c \pmod{p},\ {\rm where\ } |c| < p^{1-1/m};
$$
such $s$ exists by the Dirichlet Box Principle.  
Let $c_1,...,c_m$ be the values $c$ so produced. \footnote{Here is where we are using
the fact that $p$ is prime:  We need it to prove that such $s$ exists, and 
to extract the values of $c$ from congruences $b \equiv sc \pmod{p}$.}

Define
$$
h(n)\ =\ (\mu * f)(sn)\ =\ \SIGMA_{a+b\equiv n} \mu(sa) f(sb).
$$
We have that $h : \Z_p \to [0,1]$ and 
$$
\hat h(a)\ =\ \hat \mu(sa) \hat f(sa).
$$

Finally, define $g : \R \to \R$ to be
$$
g(\alpha)\ =\ p^{-1} \SIGMA_{1 \leq i \leq m} 
e^{-2\pi i c_i \alpha/p} \hat h(c_i),
$$
which is a truncated inverse Fourier transform of $\hat h$.  We note that if 
$|\alpha - \beta| < 1$, then since $|c_i| < p^{1-1/m}$ we deduce that 
\begin{equation} \label{g_alpha}
|g(\alpha) - g(\beta)|\ <\ p^{-1} m \Bigl |e^{2\pi i (\alpha - \beta)p^{-1/m}} - 1\Bigr | \sup_i |\hat h(c_i)| 
\ <\ \epsilon,
\end{equation}
for $p$ sufficiently large.
\bigskip

This function $g$ clearly satisfies the first property 
$$
\hat g(0)\ =\ \hat h(0)\ =\ \hat \mu(0) \hat f(0)\ =\ \hat f(0).
$$ 
(Fourier transforms are with respect to $\Z_p$).
\bigskip

Next, suppose that $n \in \Z_p$.  Then,
$$
g(n)\ =\ h(n) - p^{-1} \SIGMA_{c \neq c_1,...,c_m} 
e^{-2\pi i c n/p} \hat \mu(sc) \hat f(sc)\ =\ h(n) - \delta,
$$
where 
$$
|\delta|\ \leq\ ||\hat \mu||_1 \sup_{c \neq c_1,...,c_m} |\hat f(sc)|\ =\ 
||\hat \mu||_1 \sup_{b \in \Z_p \setminus \B'} |\hat f(b)|\ <\ \epsilon.
$$

From this, together with (\ref{g_alpha}) we have that for $\alpha \in \R$,
$g(\alpha) \in [-2\epsilon, 1 + 2\epsilon]$, as claimed by the second property
in the conclusion of the proposition.
\bigskip

Next, we observe that 
$$
\LL(g)\ =\ \LL(h) - E,\ 
$$
where 
\begin{eqnarray}
|E|\ \leq\ p^{-3} \SIGMA_{c \neq c_1,...,c_m} |\hat h(c)|^2 |\hat h(-2c)|
\ &<&\ \epsilon (\epsilon/6)^t p^{-1} ||\hat h||_2^2 \nonumber \\
&\leq&\ \epsilon^2/6. \nonumber
\end{eqnarray}
To complete the proof of the Proposition, we must relate $\LL(h)$ to $\LL(f)$:  
We begin by observing that if $b \in \B$, then 
\begin{eqnarray}
|\hat f(b) - \hat h(s^{-1}b)|\ =\ |\hat f(b)| |1 - \hat \mu(b)|\ <\ \epsilon^2 p.
\end{eqnarray}
Also, if $b \in \Z_p \setminus \B$, then
$$
|\hat f(b) - \hat h(s^{-1}b)|\ <\ 2|\hat f(b)|\ <\ 2\epsilon p.
$$
Thus,
$$
||\hat f(sa) - \hat h(a)||_\infty\ <\ 2\epsilon p.
$$

From Lemma \ref{same_fourier} with $\beta = 2\epsilon$ we conclude that 
$$
|\LL(f) - \LL(h)|\ <\ 24\epsilon.
$$
So,
$$
|\LL(f) - \LL(g)|\ <\ 25 \epsilon.\ \ \ \blacksquare
$$
\bigskip

Finally, we will require the following two technical lemmas, which are used in the proof of
Theorem \ref{second_theorem}:

\begin{lemma} \label{spread_lemma}
Suppose $p$ is prime, and suppose that $S \subseteq \Z_p$ satisfies 
$$
p/3 < |S| < 2p/5.
$$
Let $r(n)$ be the number of pairs $(s_1,s_2) \in S \times S$ such that 
$n = s_1 + s_2$.  
Then, if $T \subseteq \Z_p$, and $p$ is sufficiently large, we have 
$$
\SIGMA_{n \in T} r(n)\ <\ 0.93 |S| (|S| |T|)^{1/2}. 
$$
\end{lemma}

\noindent {\bf Proof of the Lemma.}   First, observe that if $1 \leq a \leq p-1$, then
among all subsets $S \subseteq \Z_p$ of cardinality at most $p/2$, the one which
maximizes $|\hat S(a)|$ satisfies
\begin{eqnarray}
|\hat S(a)|\ =\ \left | 1 + e^{2\pi i/p} + e^{4\pi i /p} + \cdots + 
e^{2\pi i (|S|-1)/p} \right |
\ &=&\ {|e^{2\pi i |S|/p} - 1| \over |e^{2\pi i /p} - 1|} \nonumber \\
&=&\ {|\sin(\pi |S|/p)| \over |\sin(\pi/p)|}. \nonumber
\end{eqnarray}

Since $|\theta| > \pi/3$ we have that
$$
|\sin(\theta)|\ <\ {\sin(\pi/3) |\theta| \over \pi/3}\ =\ 
{3 \sqrt{3} |\theta| \over 2\pi}.
$$
This can be seen by drawing a line passing through $(0,0)$ and 
$(\pi/3,\sin(\pi/3))$, and realizing that for $\theta > \pi/3$ we have
$\sin(\theta)$ lies below the line.
Thus, since $p/3 < |S| < 2p/5$ we deduce that for $a \neq 0$,
$$
|\hat S(a)|\ <\ {3 \sqrt{3} |S| \over 2p |\sin(\pi/p)|}\ \sim\ {3 \sqrt{3} |S| \over 2 \pi}.
$$
Thus, by Parseval,
\begin{eqnarray}
||S*S||_2^2\ =\ p^{-1} ||\hat S||_4^4\ &\leq&\ p^{-2}|S|^4+ p^{-1}(||\hat S||_2^2 - p^{-1} |S|^2) 
\sup_{a \neq 0} |\hat S(a)|^2 \nonumber \\
&<&\ 0.856 p^{-1} |S|^3, \nonumber
\end{eqnarray} 
for $p$ sufficiently large.

By Cauchy-Schwarz we have that 
\begin{eqnarray}
\SIGMA_{n \in T} r(n)\ &\leq&\ |T|^{1/2} \left ( \SIGMA_n r(n)^2 \right )^{1/2}
\nonumber \\
&=&\  |T|^{1/2} p^{1/2} ||S*S||_2 \nonumber \\
&<&\ 0.93 |S| ( |S| |T|)^{1/2}.\ \ \ \ \ \blacksquare \nonumber
\end{eqnarray}
\bigskip

\begin{lemma} \label{complement_lemma}  
Suppose $N \geq 3$ is odd, and suppose 
$A \subseteq \Z_N$, $|A| = \upsilon N$.
Let $A'$ denote the complement of $A$.  Then, 
$$
\LL(A) + \LL(A')\ =\ 3\upsilon^2 - 3\upsilon + 1
$$
\end{lemma}

\noindent {\bf Proof.}  The proof is an immediate consequence of the fact that 
$\hat A'(0) = (1-\upsilon)N$, together with $\hat A(a) = - \hat A'(a)$
for $1 \leq a \leq N-1$.  For then, we have 
\begin{eqnarray}
\LL(A) + \LL(A')\ &=&\ N^{-3} \SIGMA_a 
\hat A(a)^2 \hat A(-2a) + \hat A'(a) \hat A'(-2a) \nonumber \\
\ &=&\ \upsilon^3 + (1-\upsilon)^3 \nonumber \\
&=&\ 3\upsilon^2 - 3\upsilon + 1. \ \ \blacksquare \nonumber
\end{eqnarray}
\bigskip


\section{Proof of Theorem \ref{main_theorem}}

To prove the theorem it suffices to show that for every 
$0 < \epsilon, \upsilon < 1$, every pair of primes $p,r$ with  
$r > p^3 > p_0(\epsilon)$, and every function $f : \Z_p \to [0,1]$ satisfying
$\EE(f) \geq \upsilon$, there exists 
a function $\ell : \Z_r \to [0,1]$ satisfying $\EE(\ell) \geq \upsilon$, such that 
\begin{equation} \label{AB_target}
\LL(\ell)\ <\ \LL(f) + \epsilon
\end{equation}
This then implies
$$
\rho(\upsilon,r)\ <\ \rho(\upsilon,p) + \epsilon,
$$
and then our theorem follows (because then $\rho(r,\upsilon)$ is 
approximately decreasing as $r$ runs through the primes.)
\bigskip

To prove (\ref{AB_target}), let $f : \Z_p \to [0,1]$ satisfy $\EE(f) \geq \upsilon$.
Then, applying Proposition \ref{measure_prop} 
we deduce that there is a map $g : \R \to \R$ satisfying the
conclusion of that proposition.  Let $c_1,...,c_m$, $|c_i| < p^{1-1/m}$
be as in the proposition.

Define
$$
h(\alpha)\ =\ p^{-1} \SIGMA_{1 \leq i \leq m} e^{-2\pi i \alpha c_i/r} \hat g(c_i)\ =\ 
g(\alpha p / r)\ \in\ [-2\epsilon, 1+2\epsilon].
$$
If we restrict to integer values of $\alpha$, then we have that $h$ has the following
properties
\bigskip

$\bullet$  $h : \Z_r \to [-2\epsilon, 1+2\epsilon]$.

$\bullet$  $\EE(h)\ =\ \EE(g)\ \geq\ \upsilon r$.  (Here, $\EE(g)$ is computed by 
restricting to $g : \{0,...,p-1\} \to \R$.)

$\bullet$  For $|a| < r/2$ we have $\hat h(a) \neq 0$ if and only if $a = c_i$ for some $i$,
where $|c_i| < p^{1-1/m}$, in which case $\hat h(c_i) = r \hat g(c_i)/p$.
\bigskip

From the third conclusion we get that
$$
\LL(h)\ =\ r^{-3} \SIGMA_{1 \leq i \leq m} \hat h(c_i)^2 \hat h(-2c_i)\ =\ \LL(g).
$$ 
Then, from the final conclusion in Proposition \ref{measure_prop} we have that 
\begin{equation} \label{upper_T}
\LL(h)\ <\ \LL(f) + 25\epsilon.
\end{equation}
\bigskip

This would be the end of the proof of our theorem were it not for the fact that 
$h : \Z_r \to [-2\epsilon, 1 + 2\epsilon]$, instead of $\Z_r \to \{0,1\}$.  This is easily
fixed:  First, we let $\ell_0 : \Z_r \to [0,1]$ be defined by
$$
\ell_0(n)\ =\ \left \{ \begin{array}{rl} h(n),\ &\ {\rm if\ } h(n) \in [0,1]; \\
0,\ &\ {\rm if\ } h(n) < 0; \\
1,\ &\ {\rm if\ } h(n) > 1.\end{array}\right.
$$
We have that 
$$
|\ell_0(n) - h(n)| \leq 2\epsilon,\ {\rm and\ therefore\ } ||\hat \ell_0 - \hat h||_\infty < 2\epsilon r.
$$
It is clear that by reassigning some of the values of $\ell_0(n)$ we can produce a
map $\ell : \Z_r \to [0,1]$ such that 
\footnote{If $\hat \ell_0(0) > \hat h(0)$, then we reassign some of the $n$ where
$\ell_0(n) = 1$ to $0$, so that we then get $\hat h(0) \leq \hat \ell_0(0) < \hat h(0) + 1$, and
then we change one more $n$ where $\hat \ell_0(n) = 0$ to produce $\ell : \Z_r \to [0,1]$
satisfying $\hat \ell(0) = \hat h(0)$; likewise, if $\hat \ell_0(0) < \hat h(0)$, we reassign
some values where $\hat \ell_0(n) = 0$ to $1$.} 
$$
\EE(\ell)\ =\ \EE(h),\ {\rm and\ } ||\hat \ell - \hat h||_\infty < 4\epsilon r.
$$
From Lemma \ref{same_fourier} we then deduce
$$
|\LL(\ell) - \LL(h)|\ <\ 48 \epsilon;
$$
and so,
$$
\EE(\ell)\ =\ \EE(f),\ {\rm and\ }\LL(\ell)\ <\ \LL(f) + 73 \epsilon.
$$
Our theorem is now proved on rescaling the $73\epsilon$ to $\epsilon$. \ \ \ $\blacksquare$

\section{Proof of Theorem \ref{second_theorem}}

A consequence of Lemma \ref{complement_lemma} 
is that for a given density $\upsilon$, the
sets $A \subseteq \Z_N$ which minimize $\LL(A)$ are exactly those which 
maximize $\LL(A')$.  If $3 | N$ and $\upsilon  = 2/3$, clearly if we let 
$A'$ be the multiplies of $3$ modulo $N$, then $\LL(A')$ is maximized and
therefore $\LL(A)$ is minimized.  In this case, for every pair $m,m+d \in A'$ we have 
$m+2d \in A'$, and so $\LL(A') = (1-\upsilon)^2$.  By the above lemma,
$$
\LL(A)\ =\ 3\upsilon^2 - 3\upsilon + 1 - (1-\upsilon)^2
\ =\ 2 \upsilon^2 - \upsilon\ =\ 2/9.
$$
So, 
$$
\rho(2/3,N)\ =\ 2/9.
$$

The idea now is to show that
$$
\lim_{p \to \infty \atop p\ {\rm prime}} \rho(2/3,p)\ \neq\ 2/9.
$$

Suppose $p \equiv 1 \pmod{3}$ and that 
$A \subseteq Z_p$ minimizes $\LL(A)$ subject to $|A| = (2p+1)/3$.
Let $S = \Z_p \setminus A$, and note that $|S| = (p-1)/3$.  
Let $T = 2*S = \{2s : s\in S\}$.  

Now, if $r(n)$ is the number of pairs $(s_1,s_2) \in S \times S$
satisfying $s_1 + s_2 = n$, then by Lemma \ref{spread_lemma} we have 
$$
\LL(T)\ =\ p^{-2} \sum_{n \in T} r(n)\ <\ 0.93 p^{-2} |S| (|S| |T|)^{1/2}\ \leq\ 
0.93/ 9,
$$
for all $p$ sufficiently large.  So, by Lemma  \ref{complement_lemma} we have
that
$$
\LL(A)\ >\ 0.23,
$$
and therefore
$$
\rho(2/3,p)\ >\ 0.23\ >\ 2/9
$$
for all sufficiently large primes $p \equiv 1 \pmod{3}$.
This finishes the proof of the theorem.\ \ \ \ $\blacksquare$ 
\bigskip

\section{Acknowledgements}

I would like to thank Ben Green for the question, as well as for suggesting the
proof of Theorem 1, which was a modification of an earlier proof of the author.

\end{document}